\documentclass[preprint,12pt]{elsarticle}

\usepackage{amsmath,amssymb,amsthm}
\usepackage{tikz}
\usetikzlibrary{positioning}
\usepackage[hmargin=2.6cm,vmargin=2.6cm]{geometry}
\usepackage{hyperref}
\usepackage{setspace}
\newtheorem{theorem}{Theorem}[section]
\newtheorem{lemma}[theorem]{Lemma}
\newtheorem{proposition}[theorem]{Proposition}
\newtheorem{corollary}[theorem]{Corollary}
\newtheorem{claim}[theorem]{Claim}
\theoremstyle{definition}

\newtheorem{example}[theorem]{Example}
\newtheorem{remark}[theorem]{Remark}

\newcommand{\FF}{\mathbb{F}}
\newcommand{\RR}{\mathbb{R}}
\newcommand{\CC}{\mathbb{C}}
\newcommand{\cI}{\mathcal{I}}
\newcommand{\corank}{\mathrm{corank}}
\newcommand{\rank}{\mathrm{rank}}

\begin{document}

\begin{frontmatter}

\title{The minimum forcing number of graphs and involutory weighted adjacency matrices\tnoteref{fund}}
\tnotetext[fund]{This work is supported by NSFC (grant nos. 12271229, 12671409).}

\author[lzu]{Junlin Tian}
% \ead{...} % TODO: email of the first author
\author[lzu]{Heping Zhang\corref{cor1}}
\ead{zhanghp@lzu.edu.cn}
\cortext[cor1]{Corresponding author}
\address[lzu]{School of Mathematics and Statistics, Lanzhou University, Lanzhou,
Gansu 730000, PR China}

\begin{abstract}
Let $G$ be a graph with a perfect matching. The minimum forcing number $f(G)$ of $G$ is
the smallest size of a matching in $G$  that is contained in exactly one perfect matching. For the $d$-hypercube $Q_d$ ($d\geq 2$)
Diwan (2019) determined $f(Q_d)=2^{d-2}$ by using a linear
algebraic method, and Mohammadian (2022) generalized this method and   showed that $f(G\square K_2)=n/2$  when a bipartite graph $G$ admits order $n$ and
an involutory weighted adjacency matrix over a field $\FF$, where $\square$ denotes the Cartesian product. In this paper, we further
 prove that if $\operatorname{char}\mathbb{F}\neq 2$, then
$\frac{n}{2}\le f(G\mathbin{\square} P_k)\le \frac{n}{2}\lceil\frac{k}{2}\rceil$
for every $k\ge 2$, where $P_k$ is the $k$-vertex path; this recovers Mohammadian's result when $k=2$, and the lower bound is attained for $k=3$.
Moreover, we extend this  method to  non-bipartite graphs. Recursively constructing  the adjacency matrix of the
folded hypercube $FQ_n$ over the binary field $\mathbb F_2$, we prove that
$f(FQ_n\mathbin{\square} K_2)=2^{n-1}$ for every even $n$.  Combined with
Mohammadian's result for odd $n$, this establishes the same formula for
all $n\geq 1$. %We also show that
%the existence of involutory weighted adjacency matrices is preserved under the direct product of bipartite graphs, and
  %Finally, both combinatorial and algebraic arguments are given to show that
%$f(K_6^-\mathbin{\square} K_2)=3$, where $K_6^-$ is the graph obtained from $K_6$ by
%deleting a perfect matching.
\end{abstract}

\begin{keyword}
Perfect matching \sep Minimum forcing number \sep Folded hypercube \sep Weighted adjacency matrix \sep
Involutory matrix  \sep Orthogonal matrix
\end{keyword}

\end{frontmatter}

\section{Introduction}

All graphs considered in this paper are finite and simple. For a graph $G$, we denote
by $V(G)$ and $E(G)$ its vertex set and edge set, respectively.  An edge subset $M\subseteq E(G)$ is called a
\emph{matching} if no two edges of $M$ share a common end-vertex; a matching $M$ is
\emph{perfect} if it covers every vertex of $G$. For a perfect matching $M$ of $G$, a
subset $S\subseteq M$ is called a \emph{forcing set} of $M$ if $S$ is contained in exactly one
 perfect matching of $G$. The \emph{forcing number} of $M$, denoted by $f(G,M)$,
is the smallest size of a forcing set of $M$, and the \emph{minimum forcing number} of
$G$ is defined by $$f(G)=\min\{f(G,M)\mid M \text{ a perfect matching of }G\}.$$

The concept of forcing number of a perfect matching in a graph  corresponds to
the ``innate degree of freedom'' of a Kekul\'e structure of polycyclic aromatic
compounds  introduced by Klein and Randi\'c \cite{KR} in 1985, which plays an important role in resonance theory. Harary et al. \cite{HKZ}  later
used the term forcing in 1991. In general, computing the forcing number  of a perfect matching and the minimum
forcing number in a bipartite graph are both NP-complete problems \cite{AMM,AHM}. For surveys on
matching forcing, see \cite{CC,ZHLZ}.

Let $Q_n$ denote the $n$-dimensional hypercube, whose vertex set consists of all binary
vectors of length $n$, two vertices being adjacent if and only if the corresponding
vectors differ in exactly one coordinate. Diwan \cite{Diwan} first proposed a linear algebraic technique to obtain the following result, which  completely  settled  a conjecture due to
Pachter and Kim \cite{PK} and  Riddle \cite{Riddle}.

\begin{theorem}[Diwan \cite{Diwan}]\label{thm:diwan}
$f(Q_n)=2^{n-2}$ for any integer $n\ge 2$.
\end{theorem}

Let $\FF$ be a field. A \emph{weighted adjacency matrix} of a graph $G$ over $\FF$ is a
square matrix $A$ over $\FF$, with rows and columns indexed by $V(G)$, such that the
$(u,v)$-entry is nonzero if and only if $u$ and $v$ are adjacent in $G$. If $G$ is a
bipartite graph with a bipartition $(X,Y)$, a \emph{weighted bi-adjacency matrix} of
$G$ is a matrix $B$ over $\FF$ with rows and columns indexed by $X$ and $Y$,
respectively, such that the $(x,y)$-entry is nonzero if and only if $xy\in E(G)$.

Diwan's method is based on the following  lemma, which  relates the minimum
forcing number to the corank of a matrix. Recall that the \emph{corank} of an $n\times
n$ matrix $A$ is $\corank(A)=n-\rank(A)$.

%\begin{lemma}[Diwan \cite{Diwan}]\label{lem:diwan}
%Let $G$ be a graph of order $p$ with a perfect matching. If the largest order of an
%induced subgraph of $G$ having a unique perfect matching is $k$, then the forcing
%number of every perfect matching of $G$ is at least $(p-k)/2$.
%\end{lemma}

%\begin{proof}
%Let $M$ be any perfect matching of $G$ and $S$ a forcing set of $M$. Then $G-V(S)$ is
%an induced subgraph of $G$ with the unique perfect matching $M\setminus S$, so
%$k\ge 2|M\setminus S|$, and hence $|S|\ge (p-k)/2$.
%\end{proof}

\begin{lemma}[Diwan \cite{Diwan}, Mohammadian \cite{Moh}]\label{lem:corank-bip}
Let $G$ be a bipartite graph with a perfect matching. If $A$ (resp. $B$) is a weighted
(resp. bi-)adjacency matrix of $G$ over a field $\FF$, then
\begin{align}
f(G)&\ge \tfrac12\bigl(|V(G)|-\rank(A)\bigr)=\tfrac12\corank(A), \label{Ineq1}\\
f(G)&\ge \tfrac12\bigl(|V(G)|-2\rank(B)\bigr)=\corank(B). \label{Ineq2}
\end{align}
\end{lemma}

Diwan \cite{Diwan} constructed an integer weighted adjacency matrix $A$ of $Q_n$ such that its inverse
$A^{-1}$ is also an integer weighted adjacency matrix of $Q_n$, which enabled him to prove
Theorem \ref{thm:diwan}. Huang \cite{HH} constructed
a $\{0, 1,-1\}$-weighted adjacency matrix of $Q_d$ that is  row-orthogonal to confirm the famous
Sensitivity Conjecture in information theory.  Mohammadian \cite{Moh} generalized Diwan's idea to prisms of
more bipartite graphs.

A
square matrix $A$ is \emph{involutory} if $A^2=I$ (equivalently $A^{-1}=A$), and
\emph{orthogonal} if $AA^{\top}=I$. The \emph{order} of a graph $G$
is the number of its vertices. Let
\[
\cI_{\FF}(G):=\{A\mid A \text{ is an involutory weighted adjacency matrix of } G
\text{ over } \FF\}.
\]

\begin{theorem}[Mohammadian \cite{Moh}]\label{thm:moh}
Let $G$ be a bipartite graph of order $n$. If $\cI_{\FF}(G)\neq\emptyset$, then
$f(G\mathbin{\square}K_2)=n/2$.
\end{theorem}

Here the \emph{Cartesian product} $G\mathbin{\square}H$ of two graphs $G$ and $H$ is
the graph with vertex set $V(G)\times V(H)$ in which two vertices $(g_1,h_1)$ and
$(g_2,h_2)$ are adjacent if and only if either $g_1=g_2$ and $h_1h_2\in E(H)$, or
$h_1=h_2$ and $g_1g_2\in E(G)$. The graph $G\mathbin{\square}K_2$ is called the
\emph{prism} of $G$.

The \emph{folded hypercube} $FQ_n$ is the graph obtained from
$Q_n$ by adding, for every vertex, an edge joining it to its antipodal (complementary)
vertex, which was introduced by El-Amawy and Latifi \cite{EL}. $FQ_n$ is the Cayley graph of $\mathbb{Z}_2^n$ whose generating
set consists of the $n$ standard unit vectors together with the all-one vector. For odd $n$ the graph $FQ_n$ is bipartite, so the original method of
Mohammadian applies to this case: Alon and Zheng \cite{AZ} showed that $FQ_n$
admits a \emph{unitary signing}, i.e., a Hermitian
weighted adjacency matrix $A$ with $A^2=AA^{*}=(n+1)I$. Normalizing by $1/\sqrt{n+1}$
gives an involutory weighted adjacency matrix of $FQ_n$ over the complex field $\CC$, so
$\cI_{\CC}(FQ_n)\neq\emptyset$.  By  Theorem
\ref{thm:moh},   Mohammadian   obtained the following result.
\begin{theorem}\cite{Moh}\label{FQodd} $f(FQ_n\mathbin{\square}K_2)=2^{n-1}$ for odd $n$.
\end{theorem}

 The case of even $n$ remained open, since  $FQ_n$ is
non-bipartite and Theorem \ref{thm:moh} does not apply.  Mohammadian \cite{Moh} pointed out that it is  an interesting problem to determine $f(FQ_n\square K_2)$ for all
 even $n$.

 Very recently Shi and Zhang obtained the following result for circular graph $G\square C_{2k}$ by using Lemma \ref{lem:corank-bip}, where $C_{2k}$ is a cycle of length $2k$.

 \begin{theorem}[\cite{SZ}]\label{circular}
Let $G$ be a bipartite graph with $n$ vertices, $\mathbb{F}$ a field of characteristic not 2.
If $\mathcal I_{\mathbb{F}}(G)\neq\emptyset$ and integer $k\ge 2$, then
\[ {f(G\square C_{2k}) = n}. \]
\end{theorem}

In this paper we mainly extend the linear algebraic method for computing the minimum forcing
numbers to non-bipartite  graphs. In Section 2 we still consider bipartite graphs. We prove that the property of
admitting an involutory weighted adjacency matrix is preserved under the direct product
of bipartite graphs (Proposition \ref{thm:direct}), which yields
$f((G\times H)\mathbin{\square}K_2)=|V(G)||V(H)|/2$ (Corollary \ref{cor:direct}) when $G$ and $H$ have involutory weighted adjacency matrices. We then prove
that for a bipartite graph $G$ of order $n$ with $\cI_{\FF}(G)\neq\emptyset$ over a
field $\FF$ of characteristic not $2$,
\[
\frac{n}{2}\le f(G\mathbin{\square}P_k)\le \frac{n}{2}\Bigl\lceil\frac{k}{2}\Bigr\rceil
\quad\text{for every } k\ge 2,
\]
where $P_k$ denotes the $k$-vertex path (Theorem \ref{thm:path}). In contrast to
the circular graph $G\mathbin{\square}C_{2k}$, whose minimum forcing number equals $n$
 \cite{SZ},
$f(G\mathbin{\square}P_k)$ may depend on the parameter $k$. In Section 3 we extend the
linear algebraic method to some non-bipartite graphs and propose several types of generalizations. As applications, in particular we determine the minimum forcing number of the prism of
the folded hypercube $FQ_n$ for even $n$ (Theorem \ref{thm:fqn-even}), thus completing
the formula $f(FQ_n\mathbin{\square}K_2)=2^{n-1}$ for all positive integers $n$
(Corollary \ref{thm:fqn-all}). However the above binary field method is not suitable to non-bipartite
 graph of even regular degree. In Section 4, as a further example, we determine the minimum forcing number of the prism of a 4-regular non-bipartite
 graph by an algebraic method with minor combinatorial discussions.

\section{The minimum forcing number of bipartite graphs}

We first recall some concepts and terminology. For a graph $G$ with  a perfect matching, a subgraph $H$ of a graph $G$ is called \emph{nice} if
$G-V(H)$ has a perfect matching. For a matching $M$ of $G$, let $V(M)$ denote the set consisting of the endvertices of all edge of $M$.  For a subset \(S\subseteq V\), the induced subgraph by $S$, denoted by \(G[S]\), is the graph whose vertex set is $S$, and whose edge set consists of all edges in $E(G)$ with both endpoints in $S$. A \emph{perfect $2$-matching} (or \emph{elementary
spanning subgraph}) of $G$ is a spanning subgraph of $G$ each of whose components is
either $K_2$ or a cycle. In particular  a \emph{perfect $2$-matching} of $G$ is a 1-factor (or a perfect matching) if it contains no cycles.   For convenience, we refer to a forcing set of a perfect matching in a graph $G$ simply as a forcing set of $G$

%Throughout this section, for a bipartite graph $G$ with a bipartition $(X,Y)$, a
%weighted bi-adjacency matrix of $G$ means one with rows indexed by $X$ and columns
%indexed by $Y$, unless otherwise stated.
A part of the following  lemma was pointed out in
\cite{Moh}; we include a proof for completeness.

\begin{lemma}\label{lem:structure}
Let $G$ be a bipartite graph. If $\cI_{\FF}(G)\neq\emptyset$, then $G$ has a perfect matching, and
every $A\in\cI_{\FF}(G)$ has the form
\[
A=\begin{pmatrix} O & B \\ B^{-1} & O \end{pmatrix},
\]
where $B$ is a weighted bi-adjacency matrix of $G$ and the zero-nonzero pattern of
$B^{-1}$ coincides with that of the transpose $B^{\top}$.
\end{lemma}

\begin{proof}
For any $A\in\cI_{\FF}(G)$, $A$ is nonsingular, and there exists a nonzero term in its determinant expansion, which corresponds to a perfect $2$-matching.  Since $G$ contains no odd cycles, each perfect $2$-matching  contains a perfect matching of $G$. This also implies that $G$ has a  balanced bipartition $(X,Y)$. Let $n:=|V(G)|=2|X|$.  Rearranging the vertices according to the bipartition, we may
write $A=\begin{pmatrix} O & B \\ C & O\end{pmatrix}$, where $B$ and $C$ are square matrices of order $\frac{n}{2}$. Then
\[
A^2=\begin{pmatrix} BC & O \\ O & CB \end{pmatrix}=I,
\]
so $BC=CB=I_{\frac{n}{2}}$, which implies  $C=B^{-1}$.
Further, $B^{\top}$ and  $B^{-1}$ have the same zero-nonzero pattern.
\end{proof}

\subsection{Involutory weighted adjacency matrices}

A  square matrix $A$ over the complex field $\mathbb C$ is \emph{unitary} if $A^{*}A=AA^{*}=I$, where $A^{*}$
denotes the conjugate transpose of $A$.

\begin{proposition}\label{prop:orth}
If a bipartite graph $G$ admits an orthogonal weighted bi-adjacency matrix $B$ over a
field $\FF$, or a unitary weighted bi-adjacency matrix $C$ over
$\CC$, then $\cI_{\FF}(G)\neq\emptyset$ or $\cI_{\CC}(G)\neq\emptyset$, respectively.
\end{proposition}

\begin{proof}
If $B$ is orthogonal, then $B^{-1}=B^{\top}$, and
$A=\begin{pmatrix} O & B \\ B^{\top} & O\end{pmatrix}$
satisfies $A^2=\begin{pmatrix} BB^{\top} & O \\ O & B^{\top}B\end{pmatrix}=I$; since
$B^{\top}$ has the same zero-nonzero pattern as $B^{-1}$, we have
$A\in\cI_{\FF}(G)$. The unitary case is analogous with
$A'=\begin{pmatrix} O & C \\ C^{*} & O\end{pmatrix}$.
\end{proof}

Recall that the \emph{direct product} $G\times H$ of two graphs $G$ and $H$ is the
graph with vertex set $V(G)\times V(H)$ in which $(g_1,h_1)$ and $(g_2,h_2)$ are
adjacent if and only if $g_1g_2\in E(G)$ and $h_1h_2\in E(H)$. The \emph{tensor
product (or Kronecker product)} of an $n\times n$ matrix $X=(a_{ij})$ and an $m\times m$ matrix $Y$ is the
partitioned matrix $X\otimes Y=(a_{ij}Y)_{n\times n}$ of order $mn$.

\begin{proposition}\label{thm:direct}
Let $G$ and $H$ be bipartite graphs. If $\cI_{\FF}(G)\neq\emptyset$ and
$\cI_{\FF}(H)\neq\emptyset$, then $\cI_{\FF}(G\times H)\neq\emptyset$.
\end{proposition}

\begin{proof}
Suppose that $G$ and $H$ have orders $n$ and $m$, respectively. Take
$A(G)\in\cI_{\FF}(G)$ and $A(H)\in\cI_{\FF}(H)$; then $A(G)^2=I_n$ and $A(H)^2=I_m$.
By the definition of the tensor product,
\[
\bigl(A(G)\otimes A(H)\bigr)^2=A(G)^2\otimes A(H)^2=I_n\otimes I_m=I_{nm},
\]
so $A(G)\otimes A(H)$ is involutory. It remains to show that $A(G)\otimes A(H)$ is a
weighted adjacency matrix of $G\times H$. For $g_i,g_v\in V(G)$ and $h_j,h_w\in V(H)$,
the $((g_i,h_j),(g_v,h_w))$-entry of $A(G)\otimes A(H)$ equals
\[
A(G)_{(g_i,g_v)}\cdot A(H)_{(h_j,h_w)},
\]
which is nonzero if and only if $g_ig_v\in E(G)$ and $h_jh_w\in E(H)$, i.e., if and
only if $(g_i,h_j)$ and $(g_v,h_w)$ are adjacent in $G\times H$. Hence
$A(G)\otimes A(H)$ is an involutory weighted adjacency matrix of $G\times H$, that is,
$A(G)\otimes A(H)\in\cI_{\FF}(G\times H)\neq\emptyset$.
\end{proof}

Since the direct product of two bipartite graphs is bipartite, % (with bipartition
%$(X_G\times X_H)\cup(Y_G\times Y_H)$ and $(X_G\times Y_H)\cup(Y_G\times X_H)$),
combining Proposition \ref{thm:direct} with Theorem \ref{thm:moh} gives the following
result.

\begin{corollary}\label{cor:direct}
Let $\FF$ be a field and let $G$ and $H$ be bipartite graphs of orders $n$ and $m$,
respectively. If $\cI_{\FF}(G)\neq\emptyset$ and $\cI_{\FF}(H)\neq\emptyset$, then
\[
f\bigl((G\times H)\mathbin{\square}K_2\bigr)=\frac{nm}{2}.
\]
\end{corollary}

\begin{proposition}\label{moh}
Let $\mathbb F$ be a field with $|\mathbb F|\ge 4$. Let $G$ be a bipartite graph on $n$ vertices such that  $\mathcal{I}_{\mathbb F}(G)\ne \emptyset$. Then $\mathcal{I}_{\mathbb F}(G\square K_2)\ne \emptyset$, and
\[
f(G\square Q_d)=2^{d-2}n, d\geq 1.
\]

\end{proposition}

\subsection{The minimum forcing number of \texorpdfstring{$G\mathbin{\square}P_k$}{GxPk}}

We now consider the Cartesian product of a bipartite graph $G$ and the path $P_k$ on
$k$ vertices, and show that its minimum forcing number depends not only on the order of
$G$ but also on the parameter $k$.

\begin{theorem}\label{thm:path}
Let $\FF$ be a field of characteristic not $2$, let $G=(X,Y)$ be a bipartite graph of
order $n$. If
$\cI_{\FF}(G)\neq\emptyset$, then
\begin{equation}\label{eq:path-bounds}
\frac{n}{2}\le f(G\mathbin{\square}P_k)\le \frac{n}{2}\Bigl\lceil\frac{k}{2}\Bigr\rceil,
\qquad k\ge 2.
\end{equation}
\end{theorem}

\begin{proof}
Let $(X_i,Y_i)$, $i=1,\dots,k$, be the bipartition of the $i$-th copy $G_i$ of $G$ in
$G\mathbin{\square}P_k$ corresponding to $(X,Y)$, such that every pair of corresponding
vertices in $G_i$ and $G_{i+1}$ is adjacent. Since $\cI_{\FF}(G)\neq\emptyset$, by
Lemma \ref{lem:structure} we have $|X|=|Y|=n/2$.

We first prove the lower bound
\[
f(G\mathbin{\square}P_k)\ge \frac{n}{2}.
\]
Take $A\in\cI_{\FF}(G)$. By Lemma \ref{lem:structure}, $A=\begin{pmatrix} O & B\\
B^{-1} & O\end{pmatrix}$, where $B$ is a weighted bi-adjacency matrix of $G$ from $X$
to $Y$.  $G\mathbin{\square}P_k$ is a bipartite  graph with a bipartition
$(X_1\cup Y_2\cup X_3\cup\cdots,\ Y_1\cup X_2\cup Y_3\cup\cdots)$. We construct a
weighted bi-adjacency matrix $R$ of $G\mathbin{\square}P_k$ which is a $k\times k$
partitioned matrix $(R_{ij})_{k\times k}$, each block $R_{ij}$ being a square matrix of
order $n/2$; the block rows are indexed by $X_1,Y_2,X_3,\dots$ and the block columns by
$Y_1,X_2,Y_3,\dots$. Denote by $R^{(i)}$ the $i$-th block row of $R$, and by $I$ the
identity matrix of order $n/2$. Since $\FF$ has characteristic not $2$, the matrices
$2B$ and $2B^{-1}$ have the same zero-nonzero pattern as $B$ and $B^{-1}$,
respectively.

\medskip
%\noindent\textbf{Case 1: $k=2$.} Then $G\mathbin{\square}P_2$ is
%the prism of $G$, and $f(G\mathbin{\square}P_2)=n/2$ by Theorem \ref{thm:moh}.

%\medskip
\noindent\textbf{Case 1: $k\equiv 0\pmod 2$, $k\ge 2$.} We can write $R$ as follows.

\begin{equation}\label{R1}
   R:=\bordermatrix{%
   & Y_1    & X_2	  & Y_3		 & X_4     &	\cdots	 & X_{k-2}& Y_{k-1}& X_{k}\cr		
X_1&B & I & O & O & \cdots & O & O & O\cr
Y_2 &I & 2B^{-1} & I & O & \cdots & O & O & O\cr
X_3&O & I & 2B & I & \cdots & O & O & O\cr
Y_4& O & O & I & 2B^{-1} & \cdots & O & O & O\cr
\vdots&\vdots & \vdots & \vdots & \vdots & \ddots & \vdots & \vdots & \vdots\cr
Y_{k-2}&O & O & O & O & \cdots & 2B^{-1} & I & O\cr
X_{k-1}&O & O & O & O & \cdots & I & 2B & I\cr
Y_k&O & O & O & O & \cdots & O & I & B^{-1}
}.
\end{equation}
The diagonal blocks of $R$ satisfy $R_{ii}=2B^{-1}$ or $2B$ according to
whether $i$ is even or odd, except $R_{11}=B$ and $R_{kk}=B^{-1}$; the off-diagonal
blocks satisfy $R_{i-1,i}=R_{i,i-1}=I$ for $2\le i\le k$; all other blocks are zero.
Thus each block row (column) of $R$ has exactly three nonzero blocks, except the first
and the last which have two.
 We claim that
\begin{equation}\label{eq:even-rel}
R^{(1)}=BR^{(2)}+\sum_{i=1}^{k/2-1}\bigl(-R^{(2i+1)}+BR^{(2i+2)}\bigr).
\end{equation}
%Indeed, denote by $Z_j$ the $j$-th block of the right-hand side minus $R^{(1)}$.
For $k=2$, it is obvious that $R^{(1)}=BR^{(2)}$ and Eq. (\ref{eq:even-rel}) holds. Next let $k\geq 4$. From the right side of Eq. (\ref{eq:even-rel}) we can see that Clearly $R_{11}=B=B\cdot I$ and $R_{12}=I=B\cdot 2B^{-1}-I$.
For $3\le j\le k$ there are only three possibilities:
\begin{itemize}
\item $j$ odd: $R_{1j}=O=BI-2B+BI$;
\item $j<k$ even, $j\neq k$: $R_{1j}==O=-I+B\cdot 2B^{-1}-I$;
\item $j=k$: $R_{1k}=O=-I+BB^{-1}$.
\end{itemize}
Hence \eqref{eq:even-rel} holds, which implies that the first $n/2$ rows of $R$ are
linear combinations of the last $(kn-n)/2$ rows of $R$. Therefore
$\rank(R)\le \frac{kn-n}{2}$.

\medskip
\noindent\textbf{Case 2: $k\equiv 1\pmod 2$, $k\ge 3$.} Define $R$ as in Case 2 as follows. Note that  there are odd number of diagonal blocks of $R$ and the diagonal blocks  satisfy $R_{ii}=2B^{-1}$ or $2B$ according to
whether $i$ is even or odd, except $R_{11}=R_{kk}=B$.

\begin{equation}\label{R2}
   R:=\bordermatrix{%
           & Y_1    & X_2	  & Y_3		 & X_4     &	\cdots	 & X_{k-3}& Y_{k-2}& X_{k-1}&Y_k\cr		
X_1&B & I & O & O & \cdots & O & O & O & O\cr
Y_2&I & 2B^{-1} & I & O & \cdots & O & O & O & O\cr
X_3&O & I & 2B & I & \cdots & O & O & O & O\cr
Y_4& O& O & I & 2B^{-1} & \cdots & O & O & O & O\cr
\vdots &\vdots & \vdots & \vdots & \vdots & \ddots & \vdots & \vdots & \vdots & \vdots\cr
Y_{k-3}& O & O & O & O & \cdots & 2B^{-1} & I & O & O\cr
X_{k-2}&O & O & O & O & \cdots & I & 2B & I & O\cr
Y_{k-1}&O & O & O & O & \cdots & O & I & 2B^{-1} & I\cr
X_k&O & O & O & O & \cdots & O & O & I & B
}.
\end{equation}

\noindent We claim that
\begin{equation}\label{eq:odd-rel}
R^{(1)}=\sum_{i=1}^{(k-1)/2}\bigl(BR^{(2i)}-R^{(2i+1)}\bigr).
\end{equation}
Again, we can see that
$R_{11}=B=B\cdot I$ and $R_{12}=I=B\cdot 2B^{-1}-I$, and for $3\le j\le k$ there are  also the following  three cases.
\begin{itemize}
\item $j<k$ odd, $j\neq k$: $R_{1j}=O=BI-2B+BI$;
\item $j$ even: $R_{1j}=O=-I+B\cdot 2B^{-1}-I$;
\item $j=k$: $R_{1k}=O=B-BI$.
\end{itemize}
Hence Eq. \eqref{eq:odd-rel} holds and $\rank(R)\le\frac{kn-n}{2}$.

\medskip
In all cases $R$ is a weighted bi-adjacency matrix of $G\mathbin{\square}P_k$ over
$\FF$ (here the assumption $\mathrm{char}(\FF)\neq 2$ is essential). Combining the
above bounds with Lemma \ref{lem:corank-bip}, we obtain
\[
f(G\mathbin{\square}P_k)\ge \frac{|V(G\mathbin{\square}P_k)|}{2}-\rank(R)
\ge \frac{kn}{2}-\frac{kn-n}{2}=\frac{n}{2}.
\]

Now we prove the upper bound: $f(G\mathbin{\square}P_k)\le\frac{n}{2}\lceil\frac{k}{2}\rceil$.
If $k$ is even, let $M$ be a matching of $G\mathbin{\square}P_k$ consisting of the
edges joining a vertex of $Y_{2i-1}$ to the corresponding vertex of $Y_{2i}$ for
$i=1,\dots,k/2$. Then $|M|=\frac{n}{2}\cdot\frac{k}{2}$,
and $G\mathbin{\square}P_k-V(M)$ is the disjoint union of   \(n/2\) copies of \(P_k\), which has a unique perfect matching since \(k\) is even. Hence $M$ is uniquely
extendable to a perfect matching of $G\mathbin{\square}P_k$. If $k$ is odd, let $M'$ be
the matching consisting of the edges between $Y_{2i-1}$ and $Y_{2i}$ for
$i=1,\dots,\frac{k-1}{2}$, together with a perfect matching (between $X_k$ and $Y_k$) of the $k$-th copy $G_k$ of  $G$ by Lemma \ref{lem:structure}. Then
$|M'|=\frac{n}{2}\cdot\frac{k-1}{2}+\frac{n}{2}=\frac{n}{2}\lceil \frac{k}{2}\rceil$, and
$G\mathbin{\square}P_k-V(M')$ has a unique perfect matching. Hence $M'$ is uniquely
extendable to a perfect matching of $G\mathbin{\square}P_k$. In both cases
$f(G\mathbin{\square}P_k)\le\frac{n}{2}\lceil \frac{k}{2}\rceil$, and \eqref{eq:path-bounds}
follows.
\end{proof}

\begin{remark} When $k=2$, Theorem  \ref{thm:path} implies Theorem \ref{thm:moh} without the condition $\operatorname{char}(\mathbb F)\not=2$ since the weighted bi-adjacency matrix $R$ has no $2B$ or $2B^{-1}$ blocks. When $k\geq 3$ is odd, the upper bound can be improved:
\begin{equation}\label{improv}f(G\square P_k)\leq \frac{n}{2}\frac{k-1}{2}+f(G)\end{equation}
since it is sufficient to take a minimum forcing set in the $k$-copy $G_k$ to replace its perfect matching in constructing a forcing set $M'$ of the entire graph. Further $f(G)\leq \frac{n}{2}-1$.
\end{remark}

\begin{remark}
For $G=K_2$, Zhao and Zhang \cite{ZZ} showed that $f(K_2\square P_k)=\lceil \frac{k-1}{3}\rceil$. This shows indeed that the
minimum forcing number of $G\square P_k$ depends on the parameter $k$ as well, which is quite different from
 the circular graphs $G\mathbin{\square}C_{2k}$ (see Theorem \ref{circular}).
%minimum forcing number equals the order of $G$ whenever $\cI_{\FF}(G)\neq\emptyset$,
%i.e., it depends only on the size of $G$. Theorem \ref{thm:path} shows that the minimum
%forcing number of $G\mathbin{\square}P_k$ d, which
%is analogous to the results of Zhao and Zhang \cite{ZZ} on the minimum forcing numbers
%of grid graphs.
\end{remark}

\begin{proposition} Let $\mathbb F$ be a field with $|F|\ge 4$ and $\operatorname{char}(F) \neq 2$. Let $G$ be a bipartite graph of order $n$ such that  $\mathcal{I}_{\mathbb F}(G)\neq \emptyset$. Then
\[
f(G\square K_2\square P_3) =n.
\]
\end{proposition}
\begin{proof} By Proposition \ref{moh}, $\mathcal{I}_{\mathbb F}(G \square K_2)\neq \emptyset$.
By Theorem \ref{thm:path} we have
\[
f(G \square K_2 \square P_3) \ge   \frac{|V(G)|\times 2}{2}=n.
\].
%
%\[
%R: = \begin{pmatrix}
%A & I & 0 \\
%I & 2A^{-1} & -I \\
%0 & -I & A
%\end{pmatrix}.
%\]

%\noindent Since
%\[
%\begin{pmatrix} A & I & 0 \end{pmatrix}
%= A\begin{pmatrix} I & 2A^{-1} & -I \end{pmatrix}
%+ I\begin{pmatrix} 0 & -I & A \end{pmatrix},
%\]

Next, to prove the reverse inequality, we construct a matching $M$ of size $n$ in $G' := G\square K_2\square P_3$, uniquely determining a perfect matching of $G'$. Let $G$ admit a bipartition $(X,Y)$ such that $|X|=|Y|=\frac{n}{2}$. Obviously, $K_2\square P_3$ is a $2\times 3$ grid with the vertices labelled by $(i,j)$, $1\le i\le 2,1\le j\le 3$. Let $G_{ij}=(X_{ij},Y_{ij})$ be the copy of $G$ corresponding to $(i,j)$. Take a matching $M$ of $G'$ consisting of edges between $X_{11}$ and $X_{12}$ and edges between $Y_{12}$ and $Y_{13}$. Then $|M|=n$, and $G'-V(M)$ has a unique perfect matching consisting of all edges between $Y_{11}$ and $Y_{21}$, between $X_{21}$ and $X_{22}$, between $Y_{22}$ and $Y_{23}$ and between $X_{23}$ and $X_{13}$. So $M$ is a forcing set for $G'$ and $f(G')\leq n$.
\end{proof}

\begin{remark} The above proposition shows that the lower bound on \(f(G\square P_k)\) in (\ref{eq:path-bounds}) can be achieved not only for \(k=2\) but also for  \(k=3\). We want to know whether the general upper bound in (\ref{eq:path-bounds}) or an improvement in (\ref{improv})  can be achieved.

\end{remark}

\section{The minimum forcing number of non-bipartite graphs}

In this section we extend the linear algebraic method from bipartite graphs to
some non-bipartite graphs, and determine $f(FQ_n\mathbin{\square}K_2)$ for even $n$. However, for a non-bipartite graph $G$, the second inequality in Lemma \ref{lem:corank-bip} does not apply, and  the first one fails  for general non-bipartite graphs and  fields. We give such an example  below.

%The determinant expansion used above also shows why Lemma \ref{lem:corank-bip} does
%not extend to non-bipartite graphs, as the following example illustrates.

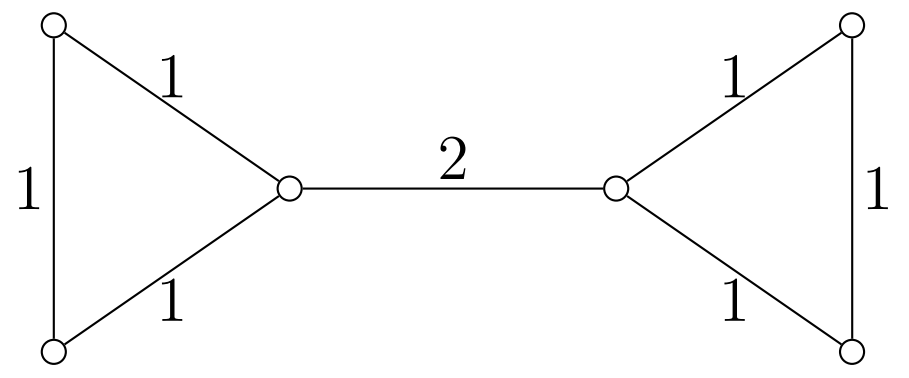
\begin{figure}[htbp]
\centering
\begin{tikzpicture}[scale=1.2, every node/.style={circle,draw,inner sep=1.6pt}]
\node (a) at (0,0.9) {};
\node (b) at (0,-0.9) {};
\node (c) at (1.3,0) {};
\node (d) at (3.1,0) {};
\node (e) at (4.4,0.9) {};
\node (f) at (4.4,-0.9) {};
\draw (a)--(b) node[midway,left,draw=none,rectangle]{$1$};
\draw (a)--(c) node[midway,above,draw=none,rectangle]{$1$};
\draw (b)--(c) node[midway,below,draw=none,rectangle]{$1$};
\draw (c)--(d) node[midway,above,draw=none,rectangle]{$2$};
\draw (d)--(e) node[midway,above,draw=none,rectangle]{$1$};
\draw (d)--(f) node[midway,below,draw=none,rectangle]{$1$};
\draw (e)--(f) node[midway,right,draw=none,rectangle]{$1$};
\end{tikzpicture}
\caption{The weighted dumbbell graph.}
\label{fig:dumbbell-wt}
\end{figure}

\begin{example}\label{ex:dumbbell}
Give the dumbbell graph $G$---two triangles joined by an edge---with the weighting
shown in Fig.~\ref{fig:dumbbell-wt}: the triangle edges have weight $1$ and the joining
edge has weight $2$. Since $G$
has a unique perfect matching, $f(G)=0$. On the other hand, $G$ has exactly
two perfect $2$-matchings: the unique perfect matching and the union of the two
triangles. Let $A(G)$ be its symmetric weighted adjacency matrix over real field by using the given edge weights.  The determinant of $A(G)$ has the following expansion.
\[
\det(A(G))=(-1)^3\cdot 2^2\cdot 1^2\cdot 1^2+(-1)^4\cdot 2^2\cdot 1^6=-4+4=0,
\]
%where in the first term a $K_2$-component whose edge has weight $w$ contributes $w^2$,
%and in the second term each triangle can be traversed in two directions.
In fact
$\rank(A(G))=5$ and $\corank(A(G))=1$, and therefore
\[
f(G)=0<\frac12=\frac12\corank(A(G)).
\]
Hence the first inequality of Lemma \ref{lem:corank-bip} does not hold for this non-bipartite
graph. \end{example}

For a general non-bipartite graph, a nice induced subgraph with a unique perfect matching may correspond to a singular weighted adjacency matrix, so the proof of Lemma \ref{lem:corank-bip} fails in such general case.

\subsection{A corank bound for general graphs}
In this subsection we  give some possible approaches that  the above linear algebra method can be extended to some special non-bipartite graphs or special fields. First we show that first lower bound on $f(G)$ in Lemma \ref{lem:corank-bip} still holds for all graphs $G$ and binary field $\mathbb F_2$, which enable us to obtain an expression for $f(FQ_n\square K_2)$.

\begin{lemma}\label{lem:det-gf2}
Let $H$ be a graph with a unique perfect matching, and let $A(H)$ be a weighted
adjacency matrix of $H$ over the binary field $\FF_2$. Then $\det(A(H))=1$.
\end{lemma}

\begin{proof}
We have  the classical expansion of the determinant of an adjacency matrix (see \cite{Biggs}):
\[
\det(A(H))=\sum_{\Lambda}(-1)^{r(\Lambda)}2^{s(\Lambda)},
\]
where $\Lambda$ runs over all perfect $2$-matchings of $H$, $r(\Lambda)$ denotes the
rank of $\Lambda$ (the number of vertices minus the number of components) and
$s(\Lambda)$ the corank of $\Lambda$ (the number of cycles of $\Lambda$). If
$s(\Lambda)\ge 1$, then $2^{s(\Lambda)}$ is even and vanishes over $\FF_2$; hence
$\det(A(H))\equiv\text{the number of perfect matchings of }H\equiv 1\pmod 2$.
\end{proof}

\begin{lemma}\label{lem:corank-gf2}
Let $G$ be a graph with a perfect matching and let $A(G)$ be a weighted adjacency
matrix of $G$ over $\FF_2$. Then
\[f(G)\ge\frac12\corank(A(G)).\]
\end{lemma}

\begin{proof}
Let $M$ be a perfect matching of $G$ with $f(G)=f(G,M)$, and let $S\subseteq M$ be a
minimum forcing set of $M$. Then $M\setminus S$ is the unique perfect matching of
$G':=G-V(S)$. The restriction $A|_{V(G')}$ is a principal submatrix of $A(G)$ on $V(G')$, and by
Lemma \ref{lem:det-gf2} it is nonsingular. Hence
\[
\corank(A(G))\le |V(G)|-2|M\setminus S|=|V(G)|-2(|M|-|S|)=2|S|,
\]
and $f(G)=|S|\ge\frac12\corank(A(G))$.
\end{proof}

The same argument yields the following general criterion.

\begin{lemma}\label{lem:corank-general}
Let $G$ be a graph with a perfect matching and let $A(G)$ be a weighted adjacency
matrix of $G$ over a field $\FF$. If for every nice induced subgraph $G'$ of $G$ having
a unique perfect matching, the restriction $A(G')=A(G)|_{V(G')}$ is nonsingular, then
\[
f(G)\ge \frac12\corank(A(G)).
\]
\end{lemma}

A sufficient condition for the nonsingularity required in Lemma
\ref{lem:corank-general} is as follows.

\begin{lemma}\label{lem:det-qi}
Let $G$ be a graph with a  perfect matching, and let $A(G)$
be a weighted adjacency matrix of $G$ over $\CC$ whose nonzero entries all belong to
$\{1,-1,\mathrm{i},-\mathrm{i}\}$. Then $f(G)\ge \frac12\corank(A(G))$.
\end{lemma}

\begin{proof}
By Lemma \ref{lem:corank-general}, it is sufficient to show that the restriction $A(G')=A(G)|_{V(G')}$ is nonsingular for any induced subgraph $G'$ with a unique perfect matching. Let $V(G')=\{v_1,\ldots,v_m\}$ and $A(G')=(a_{ij})_{m\times m}$. Consider the following expansion:
$\det(A(G'))=\sum_{\sigma}(-1)^{\tau(\sigma)}a_{1\sigma(1)}\cdots a_{m\sigma(m)}$.
A term for $\sigma$ is nonzero if and only if $v_iv_{\sigma(i)}\in E(G')$ for every $i$; the cycle
decomposition of such a permutation $\sigma$ then corresponds to a perfect $2$-matching
of $G$, where each cycle of length at least $3$ can be traversed in two directions. The
unique perfect matching of $G$ contributes exactly one nonzero term, and every perfect
$2$-matching containing $t\ge 1$ cycles contributes exactly $2^t$ nonzero terms. Hence
the total number of nonzero terms in the expansion is odd. Since each nonzero term lies
in $\{1,-1,\mathrm{i},-\mathrm{i}\}$, writing the sum as $p+q\mathrm{i}$ with
$p,q\in\mathbb{Z}$, the number of real terms is congruent to $p\pmod 2$ and the number
of imaginary terms is congruent to $q\pmod 2$. If $p+q\mathrm{i}=0$ then both numbers
would be even, contradicting that their sum is odd. Hence $\det(A(G'))\neq 0$.
\end{proof}

%\begin{lemma}\label{lem:cycle-free}
%Let $G$ be a graph with a perfect matching. Then the condition of Lemma
%\ref{lem:corank-general} holds in each of the
%following cases:
%\begin{enumerate}
%\item[(i)] $\FF=\FF_2$;
%\item[(ii)] all nonzero entries of $A(G)$ belong to $\{1,-1,\mathrm{i},-\mathrm{i}\}$;
%\item[(iii)] every nice induced subgraph of $G$ with a unique perfect matching has no
%perfect $2$-matching containing a cycle.
%\end{enumerate}
%\end{lemma}
%
%\begin{proof}
%Case (i) is Lemma \ref{lem:corank-gf2}, and case (ii) follows from Lemma
%\ref{lem:det-qi} applied to each nice induced subgraph with a unique perfect matching.
%In case (iii), for every such subgraph $H$ the determinant expansion of $A(G)|_{V(H)}$
%has exactly one nonzero term, namely the term corresponding to the unique perfect
%matching of $H$; hence $A(G)|_{V(H)}$ is nonsingular.
%\end{proof}

\begin{figure}[htbp]
    \centering
    \begin{tikzpicture}[
        scale=1.7,
        point/.style={circle,fill=black,inner sep=1.5pt},
        edge/.style={black,thick},
        rededge/.style={red,thick},
        blueedge/.style={blue,thick},
        label/.style={font=\small}
    ]

    % ===================== 左侧 K4：顺时针旋转45° =====================
    \coordinate (010) at (-2.0, 1.0);
    \coordinate (110) at (-1.0, 1.0);
    \coordinate (100) at (-1.0, 0.0);
    \coordinate (000) at (-2.0, 0.0);

    % 左 K4 内部所有边
    \draw[edge] (000) -- (010);
    \draw[edge] (000) -- (100);
    \draw[edge] (000) -- (110);
    \draw[edge] (100) -- (110);
    \draw[edge] (010) -- (110);
    \draw[line width=1.8pt, red] (010) -- (100);

    % ===================== 右侧 K4：中心对称翻转 =====================
    \coordinate (101) at (1.0, 0.0);
    \coordinate (001) at (2.0, 0.0);
    \coordinate (011) at (2.0, 1.0);
    \coordinate (111) at (1.0, 1.0);

    % 右 K4 内部所有边
    \draw[edge] (001) -- (101);
    \draw[edge] (001) -- (011);
    \draw[line width=1.8pt, red] (001) -- (111);
    \draw[edge] (101) -- (011);
    \draw[edge] (101) -- (111);
    \draw[edge] (011) -- (111);

    % ===================== 蓝色曲线连接（后两位相同） =====================
    \draw[edge] (110) -- (111);
    \draw[edge] (100) -- (101);
    \draw[edge] (000) to[bend right=25] (001);
    \draw[edge] (010) to[bend left=25] (011);

    % ===================== 实心顶点 =====================
    \foreach \p in {000,001,010,011,100,101,110,111} {
        \fill[black] (\p) circle (1.5pt);
    }

    % ===================== 标签 =====================
    \node[label, left]  at (000) {000};
    \node[label, below right] at (001) {001};
    \node[label, left] at (010) {010};
    \node[label, above right]  at (011) {011};

    \node[label, below right]  at (100) {100};
    \node[label, below left] at (101) {101};
    \node[label, above right] at (110) {110};
    \node[label, above left]  at (111) {111};

    \end{tikzpicture}
    \caption{$K_4 \Box K_2$: the two thick (red) edges form a forcing set}
    \label{FQ22}
\end{figure}
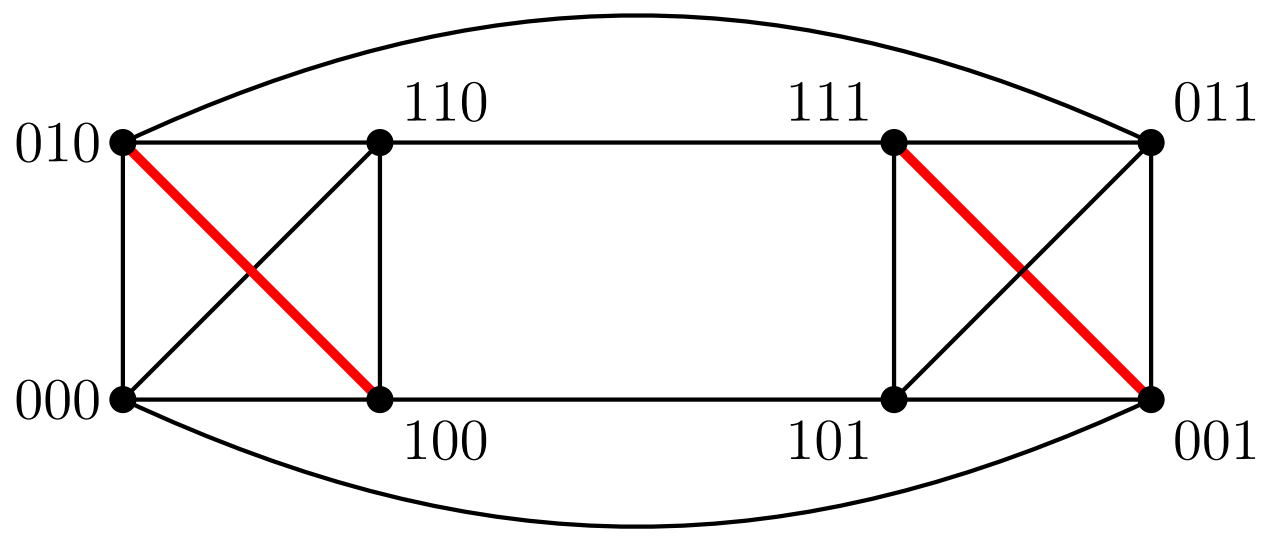

%\begin{figure}[htbp]
%\centering
%\begin{tikzpicture}[scale=1.1, every node/.style={circle,draw,inner sep=1.6pt,font=\small}]
%% copy 1: vertices 000,100,010,110
%\node (a) at (0,2) [label=above:$000$] {};
%\node (b) at (2,2) [label=above:$100$] {};
%\node (c) at (0,0) [label=below:$010$] {};
%\node (d) at (2,0) [label=below:$110$] {};
%% copy 2
%\node (a2) at (4.5,2) [label=above:$001$] {};
%\node (b2) at (6.5,2) [label=above:$101$] {};
%\node (c2) at (4.5,0) [label=below:$011$] {};
%\node (d2) at (6.5,0) [label=below:$111$] {};
%\foreach \u/\v in {a/b,a/c,a/d,b/c,b/d,c/d,a2/b2,a2/c2,a2/d2,b2/c2,b2/d2,c2/d2} \draw (\u)--(\v);
%\foreach \u/\v in {a/a2,b/b2,c/c2,d/d2} \draw (\u)--(\v);
%\draw[line width=2.2pt, red] (a)--(b);
%\draw[line width=2.2pt, red] (c2)--(d2);
%\end{tikzpicture}
%\caption{$K_4\mathbin{\square}K_2$; the two thick (red) edges form a forcing set of
%size $2$. Vertices are labelled by binary strings: the first two bits encode a vertex
%of $FQ_2=K_4$ and the last bit encodes the copy.}
%\label{fig:k4prism}
%\end{figure}

\subsection{The prism of $FQ_n$}

%Since $K_4=FQ_2$, Theorem \ref{thm:k4} solves the case $n=2$.
For  even $n$ the
graph $FQ_n$ is non-bipartite; For example, $FQ_2$ is the complete graph on the four vertices and its prism is shown in    Fig. 2. We construct the adjacency matrices of  $Q_n$ and $FQ_n$  over $\FF_2$
recursively. In what follows, $I_k$ denotes the identity matrix of order $k$, $I'_k$
denotes the matrix of order $k$ with ones on the anti-diagonal and zeros elsewhere (the
subscripts are omitted when clear from the context), and all matrices are over
$\FF_2$. Set
\[
H_1=\begin{pmatrix} 0 & 1\\ 1 & 0\end{pmatrix},\qquad
H_2=\begin{pmatrix} H_1 & I_2\\ I_2 & H_1\end{pmatrix}
=\bordermatrix{
 & 00 & 10 & 01 & 11 \cr
00 & 0 & 1 & 1 & 0 \cr
10 & 1 & 0 & 0 & 1 \cr
01 & 1 & 0 & 0 & 1 \cr
11 & 0 & 1 & 1 & 0 },
\]
\[
H'_2=\begin{pmatrix} H_1 & I_2+I'_2\\ I_2+I'_2 & H_1\end{pmatrix}
=\bordermatrix{
 & 00 & 10 & 01 & 11 \cr
00 & 0 & 1 & 1 & 1 \cr
10 & 1 & 0 & 1 & 1 \cr
01 & 1 & 1 & 0 & 1 \cr
11 & 1 & 1 & 1 & 0 }.
\]
We can see that $H_2$ is the adjacency matrix of $Q_2$ and $H'_2$ is that of
$FQ_2=K_4$.
%Similarly,
%\[
%H_3=\begin{pmatrix} H_2 & I_4\\ I_4 & H_2\end{pmatrix}
%=\bordermatrix{
% & 000 & 100 & 010 & 110 & 001 & 101 & 011 & 111 \cr
%000 & 0 & 1 & 1 & 0 & 1 & 0 & 0 & 0 \cr
%100 & 1 & 0 & 0 & 1 & 0 & 1 & 0 & 0 \cr
%010 & 1 & 0 & 0 & 1 & 0 & 0 & 1 & 0 \cr
%110 & 0 & 1 & 1 & 0 & 0 & 0 & 0 & 1 \cr
%001 & 1 & 0 & 0 & 0 & 0 & 1 & 1 & 0 \cr
%101 & 0 & 1 & 0 & 0 & 1 & 0 & 0 & 1 \cr
%011 & 0 & 0 & 1 & 0 & 1 & 0 & 0 & 1 \cr
%111 & 0 & 0 & 0 & 1 & 0 & 1 & 1 & 0 },
%\]
%\[
%H'_3=\begin{pmatrix} H_2 & I_4+I'_4\\ I_4+I'_4 & H_2\end{pmatrix}
%=\bordermatrix{
% & 000 & 100 & 010 & 110 & 001 & 101 & 011 & 111 \cr
%000 & 0 & 1 & 1 & 0 & 1 & 0 & 0 & 1 \cr
%100 & 1 & 0 & 0 & 1 & 0 & 1 & 1 & 0 \cr
%010 & 1 & 0 & 0 & 1 & 0 & 1 & 1 & 0 \cr
%110 & 0 & 1 & 1 & 0 & 1 & 0 & 0 & 1 \cr
%001 & 1 & 0 & 0 & 1 & 0 & 1 & 1 & 0 \cr
%101 & 0 & 1 & 1 & 0 & 1 & 0 & 0 & 1 \cr
%011 & 0 & 1 & 1 & 0 & 1 & 0 & 0 & 1 \cr
%111 & 1 & 0 & 0 & 1 & 0 & 1 & 1 & 0 }.
%\]
In general, $H_n$ and $H'_n$ are defined recursively by
\[
H_n=\begin{pmatrix} H_{n-1} & I\\ I & H_{n-1}\end{pmatrix},\qquad
H'_n=\begin{pmatrix} H_{n-1} & I+I'\\ I+I' & H_{n-1}\end{pmatrix}, \: n\ge 2.
\]

Order the vertices of $Q_n$ recursively as $a_10$, $a_20$, $\dots$, $a_{2^{n-1}}0$,
$a_11$, $a_21$, $\dots$, $a_{2^{n-1}}1$, where $a_1,a_2,\dots,a_{2^{n-1}}$ is the
vertex ordering of $Q_{n-1}$. For a binary string $v=v_1v_2\cdots v_n$, denote by $\bar
v=\bar v_1\bar v_2\cdots\bar v_n$ its Boolean complement.

\begin{claim}\label{claim:order}
Let $b_i$ be the binary string corresponding to the $i$-th vertex in the above ordering
of $Q_n$. Then $b_i=\bar b_{2^n+1-i}$ for $i=1,2,\dots,2^n$.
\end{claim}

\begin{proof} We use induction on $n$; the base cases $n=1,2$ are checked
directly. Next consider $n\geq 3$.
For $1\le i\le 2^{n-1}$ we have
$b_i=a_i0=\overline{a_{2^{n-1}+1-i}}0=\overline{a_{2^{n-1}+1-i}1}
=\bar b_{2^{n-1}+2^{n-1}+1-i}=\bar b_{2^n+1-i}$,
. For $2^{n-1}<j\le 2^n$,
$b_j=a_{j-2^{n-1}}1=\overline{a_{2^{n-1}+1-(j-2^{n-1})}}1
=\overline{a_{2^n+1-j}0}=\bar b_{2^n+1-j}$.
\end{proof}

\begin{lemma}\label{lem:H-construction}
Under the above vertex ordering, $H_n$ is the adjacency matrix of $Q_n$ and $H'_n$ is
the adjacency matrix of $FQ_n$ over $\FF_2$.
\end{lemma}

\begin{proof}
By induction on $n$, the block form of $H_n$ says that two vertices are adjacent if and
only if they differ in exactly one coordinate. For $FQ_n$, the additional edges join
each vertex to its complement. By Claim \ref{claim:order}, the complement of the $i$-th
vertex among $\{a_10,\dots,a_{2^{n-1}}0\}$ is the $(2^{n-1}+1-i)$-th vertex among
$\{a_11,\dots,a_{2^{n-1}}1\}$, and vice versa; these are exactly the positions of the
ones of $I'$.
% (recall that over $\FF_2$ the diagonal ones of $I$ correspond to no new
%adjacency, since $I$ and $I'$ overlap only on the diagonal, where $I'$ contributes the
%complement pairs within the off-diagonal blocks).
Hence $H'_n$ is the adjacency matrix of $FQ_n$.
\end{proof}

\begin{proposition}\label{prop:H-square}
\[
H_n^2=\begin{cases} I, & n \text{ is odd},\\ O, & n \text{ is even}.\end{cases}
\]
\end{proposition}

\begin{proof}
We proceed by induction on $n$. For $n=1$, $H_1^2=I$; for $n=2$, a direct computation gives
$H_2^2=\begin{pmatrix} 2I & 2I'\\ 2I' & 2I\end{pmatrix}=O$ over $\FF_2$. Assume the
statement for $n-1$ for $n\geq 3$. Then
\[
H_n^2=\begin{pmatrix} H_{n-1} & I\\ I & H_{n-1}\end{pmatrix}^2
=\begin{pmatrix} H_{n-1}^2+I & 2H_{n-1}\\ 2H_{n-1} & H_{n-1}^2+I\end{pmatrix}
=\begin{cases} I, & n \text{ is odd},\\ O, & n \text{ is even},\end{cases}
\]
since $2H_{n-1}=O$ over $\FF_2$ and, by the induction hypothesis, $H_{n-1}^2+I=I$ or
$O$ according to whether $n$ is odd or even.
\end{proof}

\begin{proposition}\label{prop:Hp-square}
\[
(H'_n)^2=\begin{cases} I, & n\ge 2 \text{ is even},\\ O, & n\ge 3 \text{ is odd}.\end{cases}
\]
\end{proposition}

\begin{proof}
We first show that
\begin{equation}\label{eq:commute}
H_{n-1}I'=I'H_{n-1}.
\end{equation}
For $n=2$ this is checked directly: $H_2I'=\begin{pmatrix} I' & I\\ I & I'\end{pmatrix}
=I'H_2$. Assume \eqref{eq:commute} holds for $n-1$ for $n\geq 3$. Since
$I'_{2^{n}}=\begin{pmatrix} O & I'_{2^{n-1}}\\ I'_{2^{n-1}} & O\end{pmatrix}$, we have
\[
H_nI'_{2^{n}}=\begin{pmatrix} H_{n-1} & I\\ I & H_{n-1}\end{pmatrix}
\begin{pmatrix} O & I'\\ I' & O\end{pmatrix}
=\begin{pmatrix} I' & H_{n-1}I'\\ H_{n-1}I' & I'\end{pmatrix},
\]
\[
I'_{2^{n}}H_n=\begin{pmatrix} O & I'\\ I' & O\end{pmatrix}
\begin{pmatrix} H_{n-1} & I\\ I & H_{n-1}\end{pmatrix}
=\begin{pmatrix} I' & I'H_{n-1}\\ I'H_{n-1} & I'\end{pmatrix},
\]
and the two are equal by the induction hypothesis.

Now we prove the proposition by induction on $n$. For $n=2$,  a direct computation gives
\[(H'_2)^2=\begin{pmatrix} I_2+2J_2 & 2J_2\\ 2J_2 & I_2+2J_2\end{pmatrix}=I.\]

Assume the
statement for $n-1$ for $n\geq 3$. We have
\[
(H'_n)^2=\begin{pmatrix} H_{n-1} & I+I'\\ I+I' & H_{n-1}\end{pmatrix}^2
=\begin{pmatrix} H_{n-1}^2+(I+I')^2 & H_{n-1}(I+I')+(I+I')H_{n-1}\\ (I+I')H_{n-1}+H_{n-1}(I+I') & H_{n-1}^2+(I+I')^2\end{pmatrix},
\]
where the off-diagonal blocks equal
$2H_{n-1}+H_{n-1}I'+I'H_{n-1}=O$ by \eqref{eq:commute}, and
$(I+I')^2=I+2I'+(I')^2=O$ since $(I')^2=I$. Hence
$(H'_n)^2=\begin{pmatrix} H_{n-1}^2 & O\\ O & H_{n-1}^2\end{pmatrix}$,
which equals $I$ for even $n\ge 2$ and $O$ for odd $n\ge 3$ by Proposition
\ref{prop:H-square}.
\end{proof}
\begin{remark} By Proposition \ref{prop:H-square}, $Q_n$ has an involutory adjacency matrix $H_n$
over the binary field $\mathbb F_2$ for all odd $n\geq 3$. So in this case we  reprove $f(Q_n)=2^{n-2}$  differently  from Diwan and Mohammadian's approaches.
%\begin{remark}
%Propositions \ref{prop:H-square} and \ref{prop:Hp-square} reflect the following
%algebraic facts: over $\FF_2$ the adjacency matrix of $Q_n$ is involutory exactly when
%$n$ is odd, and that of $FQ_n$ is involutory exactly when $n$ is even. Note that for
%even $n$ every vertex of $FQ_n$ has odd degree $n+1$, which is consistent with the
%involutoriness; for odd $n$ every vertex has even degree, and the analogous squaring
%cannot produce the identity.
%\end{remark}

\begin{theorem}\label{thm:fqn-even}
$f(FQ_n\mathbin{\square}K_2)=2^{n-1}$ for every even integer $n\ge 2$.
\end{theorem}

\begin{proof}
Let $G=FQ_n\mathbin{\square}K_2$ and denote the two copies of $FQ_n$ in $G$ by $H$ and
$H'$. We first show that $f(G)\le 2^{n-1}$. Let $M$ consist of the antipodal edges
(joining each vertex to its complement) of $H$ whose ends contain an odd number of
ones, together with the antipodal edges of $H'$ whose ends contain an even number of
ones. Since $n$ is even, the two ends of an antipodal edge have the same parity, so $M$
is a matching with $|M|=2^{n-2}+2^{n-2}=2^{n-1}$. Every vertex of $H$ outside $V(M)$
has an even number of ones; all its neighbors in $H$ (obtained by flipping one bit or
taking the complement) have an odd number of ones and hence lie in $V(M)$, except its
antipode, which is the only remaining neighbor. The same argument applies to $H'$, and
no edge between the two copies joins two remaining vertices, since corresponding
vertices in $H$ and $H'$ should have the same parity. Therefore $G-V(M)$ is a 1-factor. That is,  $M$ is uniquely
extendable to a perfect matching of $G$, and $f(G)\le 2^{n-1}$.

Next we show that $f(G)\ge 2^{n-1}$. By Lemma \ref{lem:H-construction}, $G$ has the
following adjacency matrix over $\FF_2$ (with the vertices of the two copies ordered as
in Lemma \ref{lem:H-construction}):
\[
\begin{pmatrix} H'_n & I\\ I & H'_n\end{pmatrix}
=\begin{pmatrix} H'_n\\ I\end{pmatrix}\begin{pmatrix} I & H'_n\end{pmatrix},
\]
where the factorization uses $(H'_n)^2=I$ from Proposition \ref{prop:Hp-square}. Hence
this matrix has rank at most $2^n$ and  corank at least $2^n$. By Lemma
\ref{lem:corank-gf2},
\[
f(FQ_n\mathbin{\square}K_2)\ge \frac12\corank\begin{pmatrix} H'_n & I\\ I & H'_n
\end{pmatrix}\ge 2^{n-1}.
\]
Hence $f(FQ_n\mathbin{\square}K_2)=2^{n-1}$.
\end{proof}

Combining Theorems \ref{thm:fqn-even} and \ref{FQodd}, we obtain the general formula as follows.

\begin{corollary}\label{thm:fqn-all}
$f(FQ_n\mathbin{\square}K_2)=2^{n-1}$ for every positive integer $n$.
\end{corollary}

\end{remark}

\section{Further examples for using Lemma \ref{lem:corank-general}}

There are further examples of non-bipartite graphs whose minimum forcing number can be computed by Lemma \ref{lem:corank-general}. The following structural property of graphs with a unique perfect matching will be used in
such discussions.

\begin{lemma}[Lov\'asz and Plummer \cite{LP}]\label{lem:bridge}
Let $G$ be a connected graph with a unique perfect matching $M$. Then $G$ has a bridge
 contained in the perfect matching.
\end{lemma}

\subsection{The prism  $K_4\square K_2$}

 Note that $K_4=FQ_2$. So $f(K_4\square K_2)=2$ has been already obtained in Theorem \ref{thm:fqn-even} by using the binary  field. However, for this graph any field is allowed.

%\begin{theorem}\label{thm:k4}
%$f(K_4\mathbin{\square}K_2)=2$.
%\end{theorem}
%
%\begin{proof}
%Let $G=K_4\mathbin{\square}K_2$ with the two copies of $K_4$ denoted by $K_4$ and
%$K_4'$. The two thick edges shown in Fig.~\ref{fig:k4prism} form a forcing set of $G$
%(the reader may verify that they uniquely determine a perfect matching), so
%$f(G)\le 2$.

For the lower bound, by a combinatorial argument we can show that  any  nice induced subgraph of $K_4\square K_2$ with a
unique perfect matching $M$ has no other  perfect $2$-matching except $M$. So the condition of Lemma \ref{lem:corank-general} holds for any weighted adjacency matrix of $K_4\square K_2$ over any field $\mathbb F$.
%By Lemma \ref{lem:bridge}, $H$ has a bridge $xy\in M$
%whose removal leaves two odd components $H_x\ni x$ and $H_y\ni y$. Then $H$ must have a
%pendant edge: otherwise $H_x$ and $H_y$ are both isomorphic to $K_3$, contained in the
%two different copies of $K_4$, respectively, and then at least two pairs of
%corresponding vertices of $H_x$ and $H_y$ are adjacent in $G$, contradicting that $xy$
%is a bridge of the induced subgraph $H$. Deleting the ends of a pendant edge of $H$
%yields a graph $H'$ which again has a pendant edge; repeating this argument shows that
%$H$ has $M$ as its only perfect $2$-matching, i.e., $H$ has no perfect $2$-matching
%containing a cycle.
Now take the following weighted adjacency matrix of $K_4$ over $\mathbb R$:
\[
A_1=\begin{pmatrix}
0 & 1 & 1 & 1\\
1 & 0 & -1 & 1\\
-1 & -1 & 0 & 1\\
1 & -1 & 1 & 0
\end{pmatrix}.
\]
One verifies directly that the rows of $A_1$ are pairwise orthogonal and
$A_1A_1^{\top}=3I$; hence $B:=\frac{1}{\sqrt3}A_1$ is an orthogonal weighted adjacency
matrix of $K_4$ over $\RR$. Construct the following weighted adjacency matrix of $G$:
\[
A(G)=\begin{pmatrix} B & I \\ I & B^{-1}\end{pmatrix}
=\begin{pmatrix} B \\ I\end{pmatrix}\begin{pmatrix} I & B^{-1}\end{pmatrix}.
\]
Thus $\corank(A(G))\ge 4$. By Lemma
 \ref{lem:corank-general},
$f(G)\ge\frac12\corank(A(G))\ge 2$.
%\end{proof}

\subsection{Prism $K_6^-\square K_2$}

Propositions \ref{prop:H-square} and \ref{prop:Hp-square} reflect the following
algebraic facts: $Q_n$ for odd $n$ admits an involutory  adjacency matrix over  $\FF_2$,
and $FQ_n$ for even  $n$ admits an involutory  adjacency matrix over  $\FF_2$. In fact such two types of graphs have odd regular degree.
When a graph has a vertex of even degree, the square of the adjacency matrix over $\FF_2$ can never be the identity, so  that
approach fails.

In this subsection we consider a typical  example of an even-degree regular graph. Let $K_6^-$ be the
graph obtained from the complete graph $K_6$ by deleting a perfect matching. Equivalently,
$V(K_6^-)=\{0,1,2,3,4,5\}$ and two vertices $i,j$ are adjacent if and only if
$i-j\equiv \pm 1, \pm 2\pmod 6$.  So  $K_6^-$ is non-bipartite and
$4$-regular. We first discuss its structural property and  determine $f(K_6^-\mathbin{\square}K_2)$ by an
algebraic method.

%\subsection{The combinatorial method}

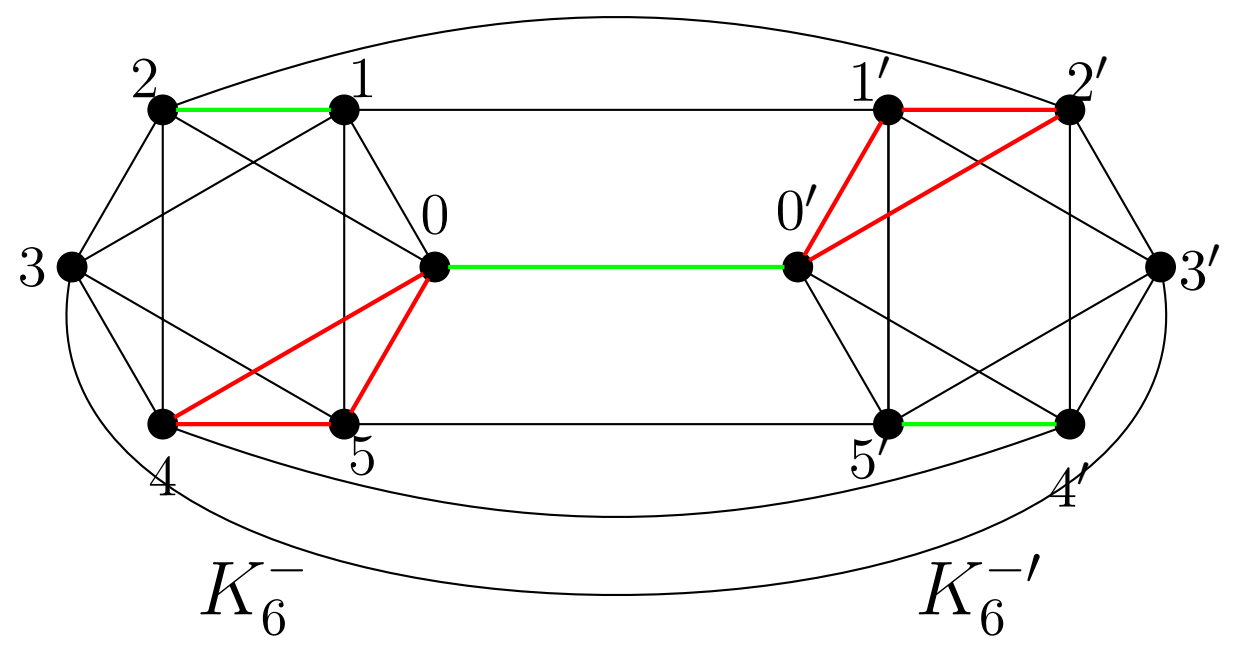
\begin{figure}[htbp]
\centering
\begin{tikzpicture}[
scale=1.2,
vertex/.style={circle,draw,thick,fill,minimum size=5pt,inner sep=0pt,outer sep=0pt},
vlab/.style={circle,draw=none,fill=none,inner sep=1pt,outer sep=0pt,font=\small,text=black}
]
% ==================== 左边：第一个 G6 ====================
\begin{scope}[xshift=-2cm]
\node[vertex] (0a) at (0:1)   {};
\node[vertex] (1a) at (60:1)  {};
\node[vertex] (2a) at (120:1) {};
\node[vertex] (3a) at (180:1) {};
\node[vertex] (4a) at (240:1) {};
\node[vertex] (5a) at (300:1) {};
\draw (2a)--(3a) (3a)--(4a)   (1a)--(0a);
\draw (1a)--(3a) (2a)--(4a) (3a)--(5a)  (5a)--(1a) (2a)--(0a);
\draw [red, thick] (0a)--(4a) ;
\draw [red, thick] (4a)--(5a);
\draw [red, thick] (5a)--(0a);
\draw [green, thick]  (2a)--(1a);
% ----- 顶点标号 -----
\node[vlab, above=4pt] at (0a) {0};
\node[vlab] at (60:1.20)   {1};
\node[vlab] at (120:1.20)  {2};
\node[vlab] at (180:1.22)  {3};
\node[vlab, below=4pt] at (4a) {4};
\node[vlab] at (300:1.20)  {5};
% 左边标签
\node[draw=none, fill=none, text=black, font=\large] at (0, -1.8) {$K_6^-$};
\end{scope}
% ==================== 右边：第二个 G6 ====================
\begin{scope}[xshift=2cm]
\node[vertex] (3b) at (0:1)   {};
\node[vertex] (2b) at (60:1)  {};
\node[vertex] (1b) at (120:1) {};
\node[vertex] (0b) at (180:1) {};
\node[vertex] (5b) at (240:1) {};
\node[vertex] (4b) at (300:1) {};
\draw (1b)--(5b)  (3b)--(4b) (0b)--(4b);
\draw (2b)--(4b)  ;
\draw (0b)--(5b)   (2b)--(3b) (3b)--(5b) (5b)--(1b) (0b)--(2b) (1b)--(3b);
\draw [red, thick] (0b)--(1b) (2b)--(1b) (0b)--(2b);
\draw [green, thick] (5b)--(4b);
% ----- 顶点标号 -----
\node[vlab, above=4pt] at (0b) {$0'$};
\node[vlab] at (120:1.20)  {$1'$};
\node[vlab] at (60:1.20)   {$2'$};
\node[vlab] at (0:1.22)    {$3'$};
\node[vlab, below=4pt] at (4b) {$4'$};
\node[vlab] at (240:1.20)  {$5'$};
% 右边标签
\node[draw=none, fill=none, text=black, font=\large] at (0, -1.8) {$K_6^{-\prime}$};
\end{scope}
% ==================== 笛卡尔积连线 ====================
\draw [green, thick] (0a) -- (0b);
\draw (1a) -- (1b);
\draw (2a) to[out=20, in=160] (2b);
\draw (3a) to[out=-100, in=-80] (3b);
\draw (4a) to[out=-20, in=-160] (4b);
\draw (5a) -- (5b);
\end{tikzpicture}
\caption{$K_6^-\mathbin{\square}K_2$; the dumbbell  induced by
$\{4,5,0,0',1',2'\}$ is  the triangles $\{0,4,5\}$ and $\{0',1',2'\}$  joined by the edge $00'$; $\{12,\ 00',\ 4'5'\}$ is a forcing set.}
\label{fig:dumbbell}
\end{figure}
\begin{lemma}\label{lem:dumbbell}
Let $G=K_6^-\mathbin{\square}K_2$ (see Fig. \ref{fig:dumbbell}). Every nice induced subgraph $H$ of $G$ having a
unique perfect matching $M$ satisfies one of the following:
\begin{enumerate}
\item[(i)]  $M$ is the only perfect $2$-matching of $H$;
\item[(ii)] $H$ is a \emph{dumbbell} such that the two  triangles  lie in the two different copies of $K_6^-$ (see Fig.~\ref{fig:dumbbell}).
\end{enumerate}
\end{lemma}

\begin{proof}
Since $K_6^-\mathbin{\square}K_2$ is $4$-connected, any induced subgraph of $G$ on at
least ten vertices is $2$-connected. By Lemma \ref{lem:bridge}, $H$ has a bridge
in $M$ and $|V(H)|\leq 8$.

Suppose that $H$ has a perfect $2$-matching containing a cycle. Since $H$ has a unique
perfect matching, every cycle in such a perfect $2$-matching is odd (an even cycle
would yield two distinct perfect matchings of $H$); since $|V(H)|$ is even, such a perfect $2$-matching has exactly two odd cycles, each of length $3$ or
$5$. Hence $6\leq |V(H)|\leq 8$

If $|V(H)|=6$,    $H$ is just a dumbbell.
   Since $H$ is induced, the two
triangles lie in the two different copies of $K_6^-$, respectively (two disjoint
triangles in one copy of $K_6^-$ would be joined by further edges), and they share
exactly one pair of corresponding vertices, as shown in Fig.~\ref{fig:dumbbell}.

If $|V(H)|=8$, then   take a bridge $uv$ of $H$ in $M$ such that $H-uv$ is the union of two disjoint  odd subgraphs $H_u$ of order 3 and $H_v$ of order 5 such that $u\in V(H_u)$ and $v\in V(H_v)$.
%$H$ has a pendant edge $xy$, then $xy\in M$ and $xy$ forms a $K_2$-component of  the
%perfect $2$-matching. Then $H':=H-x-y$ has a unique perfect matching and  two odd
%cycles of $H$. Hence $H'$
%is just a dumbbell. Since $H'$ is induced, the two
%triangles lie in the two different copies of $K_6^-$, respectively (two disjoint
%triangles in one copy of $K_6^-$ would be joined by further edges), and they share
%exactly one pair of corresponding vertices, as shown in Fig.~\ref{fig:dumbbell}. But
%each of the six vertices of $G$ outside $V(H')$ has at leat two neighbors in $V(H')$, contradicting that
%$xy$ is a pendant edge in $H$.
%Hence $\delta(H)\ge 2$ and $H$ is connected. $H-uv$ has exactly two odd components $H_u$ and $H_v$ such that $u\in V(H_u)$ and $v\in V(H_v)$;
%If $|V(H)|=8$, then one of $H_u,H_v$ has three vertices
%and the other has five; say $|H_u|=3$ and $|H_v|=5$. Since $\delta(H)\ge 2$, $H_u$ is a
%triangle, and
Since the bridge $uv$ is the only edge between $H_u$ and $H_v$, $H_u$ and $H_v$ lie in different copies of $K_6^-$. But there are at least two pairs of
corresponding vertices of $H_u$ and $H_v$, which implies that there are  at least two edges between $H_u$ and $H_v$, a contradiction.
\end{proof}

\begin{theorem}\label{thm:g6}
$f(K_6^-\mathbin{\square}K_2)=3$.
\end{theorem}
\begin{proof}
%Let $G=K_6^-\mathbin{\square}K_2$. Let $M$ be a perfect matching of $G$ with
%$f(G)=f(G,M)$ and let $S$ be a minimum forcing set of $M$, so that $H:=G-V(S)$ has the
%unique perfect matching $M\setminus S$. By Lemma \ref{lem:dumbbell}, either $H$ has no
%perfect $2$-matching containing a cycle, in which case $|V(H)|\le 6$ by Corollary
%\ref{cor:six}, or $H$ is a dumbbell and $|V(H)|=6$. In both cases $|V(S)|\ge 12-6=6$,
%i.e., $f(G)=|S|\ge 3$.
Let $M'=\{12,\ 00',\ 4'5'\}$ (see Fig.~\ref{fig:dumbbell}). Then
$G-V(M')$ is the dumbbell induced by $\{3,4,5,1',2',3'\}$, whose unique perfect
matching is $\{45,33',1'2'\}$. Hence $M'$ is uniquely
extendable to a perfect matching of $G$, and
$f(G)\le |M'|=3$. %Therefore $f(K_6^-\mathbin{\square}K_2)=3$.
%\end{proof}

%\begin{figure}[htbp]
%\centering
%\begin{tikzpicture}[scale=1.05, every node/.style={circle,draw,inner sep=1.6pt,font=\small}]
%\foreach \i/\ang in {0/90, 1/30, 2/-30, 3/-90, 4/-150, 5/150} {
%  \node (v\i) at (\ang:1.6) [label={\ang}:$\i$] {};
%}
%\foreach \i/\j in {0/1,1/2,2/3,3/4,4/5,5/0,0/2,1/3,2/4,3/5,4/0,5/1} \draw (v\i)--(v\j);
%\begin{scope}[xshift=5.2cm]
%\foreach \i/\ang in {0/90, 1/30, 2/-30, 3/-90, 4/-150, 5/150} {
%  \node (w\i) at (\ang:1.6) [label={\ang}:$\i'$] {};
%}
%\foreach \i/\j in {0/1,1/2,2/3,3/4,4/5,5/0,0/2,1/3,2/4,3/5,4/0,5/1} \draw (w\i)--(w\j);
%\end{scope}
%\foreach \i in {0,...,5} \draw (v\i)--(w\i);
%% forcing set M' = {12, 33', 4'5'} in green
%\draw[line width=2.2pt, green!60!black] (v1)--(v2);
%\draw[line width=2.2pt, green!60!black] (v3)--(w3);
%\draw[line width=2.2pt, green!60!black] (w4)--(w5);
%\end{tikzpicture}
%\caption{A forcing set $M'=\{12,33',4'5'\}$ (thick green edges) of size $3$ in
%$K_6^-\mathbin{\square}K_2$: deleting its ends leaves the dumbbell of
%Fig.~\ref{fig:dumbbell}, which has a unique perfect matching.}
%\label{fig:g6forcing}
%\end{figure}

%\subsection{The algebraic method}

We now prove the lower bound of Theorem \ref{thm:g6} by the linear algebraic method.
Construct the following weighted adjacency matrix of $K_6^-$ over $\mathbb R$:
\[
A(K_6^-)=\frac{1}{2}\begin{pmatrix}
0 & 1 & 1 & 0 & 1 & 1\\
1 & 0 & -1 & 1 & 0 & 1\\
-1 & -1 & 0 & 1 & 1 & 0\\
0 & 1 & -1 & 0 & 1 & -1\\
1 & 0 & 1 & 1 & 0 & -1\\
1 & -1 & 0 & -1 & 1 & 0
\end{pmatrix}.
\]
A direct computation shows that
$A(K_6^-)A(K_6^-)^{\top}=I$. That is, $A_6:=A(K_6^-)$
is an orthogonal  matrix of $K_6^-$ over  $\RR$. We construct the following weighted adjacency matrix of $G=K_6^-\mathbin{\square}K_2$ below.
\[
A(G)=\begin{pmatrix} A_6 & I\\ I & A_6^{\top}\end{pmatrix}
=\begin{pmatrix} A_6\\ I\end{pmatrix}\begin{pmatrix} I & A_6^{\top}\end{pmatrix},
\]
which implies   $\corank(A(G))\ge 6$.

Next we confirm that the condition of Lemma \ref{lem:corank-general} holds for $A(G)$.
Let $H$ be any  nice induced subgraph of $G$ with a unique
perfect matching. If $H$ has no perfect $2$-matching containing a cycle, then the determinant expansion of $A(G)|_{V(H)}$ has exactly one
nonzero term, so $A(G)|_{V(H)}$ is nonsingular. Otherwise, by Lemma \ref{lem:dumbbell} $H$ is a
dumbbell. Up to the symmetry
interchanging the two copies and the cyclic symmetry $i\mapsto i+1$ of $K_6^-$, the
dumbbell $H$ has vertex set
$\{a+1,a+2,a+3\}\cup\{(a+3)',(a+4)',(a+5)'\}$ for some $a\in\{0,1,\dots,5\}$ (indices
modulo $6$). We carry
out the computation for the first case in detail; the others are analogous.

Let $H$ be induced by $\{4,5,0,0',1',2'\}$. Then $H$ has exactly two perfect
$2$-matchings. One is its unique perfect matching $\{45,00',1'2'\}$, whose
contribution to $\det(A(G)|_{V(H)})$ is
\[
(-1)^3\cdot\Bigl(-\frac12\Bigr)\cdot\frac12\cdot 1\cdot 1\cdot
\Bigl(-\frac12\Bigr)\cdot\Bigl(-\frac12\Bigr)=\Bigl(\frac12\Bigr)^4.
\]
The other is the disjoint union of the two triangles; each triangle can be traversed in
two directions, so it contributes $2^2=4$ terms:
\[
\begin{aligned}
&(-1)^4\Bigl[\Bigl(-\tfrac12\Bigr)\tfrac12\tfrac12\tfrac12\tfrac12\Bigl(-\tfrac12\Bigr)
+\Bigl(-\tfrac12\Bigr)\tfrac12\tfrac12\Bigl(-\tfrac12\Bigr)\Bigl(-\tfrac12\Bigr)\tfrac12
+\Bigl(-\tfrac12\Bigr)\tfrac12\tfrac12\tfrac12\tfrac12\tfrac12\\
&\qquad\qquad
+\Bigl(-\tfrac12\Bigr)\Bigl(-\tfrac12\Bigr)\tfrac12\tfrac12\tfrac12\tfrac12\Bigr]
=\Bigl(\tfrac12\Bigr)^6-\Bigl(\tfrac12\Bigr)^6-\Bigl(\tfrac12\Bigr)^6
+\Bigl(\tfrac12\Bigr)^6=0.
\end{aligned}
\]
Hence $\det(A(G)|_{V(H)})=(1/2)^4+0=1/16\neq 0$. %The results for all six cases are
%listed in Table \ref{tab:dumbbell};
In each  of the other five cases the contribution of the two-triangle
perfect $2$-matching is $0$ and the determinant equals $\pm 1/16$.

%\begin{table}[htbp]
%\centering
%\caption{The six dumbbells and the determinants of the corresponding principal
%submatrices of $A(G)$.}
%\label{tab:dumbbell}
%\begin{tabular}{ccl}
%\hline
%Case & Vertex set of $H$ & $\det(A(G)|_{V(H)})$\\
%\hline
%1 & $\{4,5,0,0',1',2'\}$ & $\phantom{-}1/16$\\
%2 & $\{0,5,1,1',2',3'\}$ & $\phantom{-}1/16$\\
%3 & $\{0,1,2,2',3',4'\}$ & $-1/16$\\
%4 & $\{1,2,3,3',4',5'\}$ & $\phantom{-}1/16$\\
%5 & $\{2,3,4,4',0',5'\}$ & $\phantom{-}1/16$\\
%6 & $\{3,4,5,5',0',1'\}$ & $-1/16$\\
%\hline
%\end{tabular}
%\end{table}

Therefore $A(G)|_{V(H)}$ is nonsingular for every nice induced subgraph $H$ of $G$
with a unique perfect matching. By Lemma \ref{lem:corank-general},
$f(G)\ge \frac12\corank(A(G))$ $\ge 3$ and the theorem holds.
%Since the forcing set $M'$ of size $3$ in Fig.~\ref{fig:g6forcing} gives $f(G)\le 3$,
%we obtain again
%\[
%f(K_6^-\mathbin{\square}K_2)=3.
%\]
\end{proof}
\section*{Declaration of competing interest}
The authors declare that they have no known competing financial interests or personal
relationships that could have appeared to influence the work reported in this paper.

%\section*{Acknowledgements}
%This work is supported by NSFC (grant nos. 12271229, 12671409).

\section*{Data availability}
No data was used for the research described in the article.

%\section*{Declaration of  AI usage}
%The authors did not use AI tools to conduct the research. AI was used only to find  typos and improve some phrasing of the manuscript. The authors have reviewed and edited the content as needed and take full responsibility for the publication.

\end{document}